\documentclass[12pt]{article}

\usepackage{amssymb,latexsym}
\newtheorem{theorem}{Theorem}[section]
\newtheorem{proposition}[theorem]{Proposition}
\newtheorem{corollary}[theorem]{Corollary}
\newtheorem{definition}[theorem]{Definition}

\newenvironment{proof}{\medskip\noindent{\it Proof.\ }}{\mbox{$\Box$}\medskip}

\begin{document}

\def\cP{{\cal P}}
\def\cQ{{\cal Q}}

\def\sof{\hfill\rule{2mm}{2mm}}
\def\ls{\leq}
\def\gs{\geq}
\def\SS{\mathcal S}
\def\qq{{\bold q}}
\def\txx{{\frac1{2\sqrt{x}}}}
\def\Bn{\mathcal{P}_n}

\title{$132$-avoiding Two-stack Sortable Permutations, Fibonacci Numbers, and Pell Numbers\footnote{MR Subject Classification:  05A15}}

\author{Eric S. Egge \\
Department of Mathematics \\
Gettysburg College\\
Gettysburg, PA  17325  USA \\[4pt]
eggee@member.ams.org \\
\\
Toufik Mansour\\
LaBRI, Universit\'e Bordeaux 1,
              351 cours de la Lib\'eration\\
              33405 Talence Cedex, France\\[4pt]
toufik@labri.fr}
\maketitle

\begin{abstract}
In \cite{W20} West conjectured that there are $2(3n)!/((n+1)!(2n+1)!)$ two-stack sortable permutations on $n$ letters.
This conjecture was proved analytically by Zeilberger in \cite{Z24}.
Later, Dulucq, Gire, and Guibert \cite{DG1} gave a
combinatorial proof of this conjecture.
In the present paper we study generating functions for the number
of two-stack sortable permutations on $n$ letters avoiding (or
containing exactly once) $132$ and avoiding (or containing exactly
once) an arbitrary permutation $\tau$ on $k$ letters.
In several interesting cases this generating function can be expressed in terms of the generating function for the Fibonacci numbers or the generating function for the Pell numbers.

\medskip

\end{abstract}

\section{Introduction and Notation}

Let $S_n$ denote the set of permutations of $\{1, \ldots, n\}$, written in one-line notation, and suppose $\pi \in S_n$ and $\sigma \in S_k$.
We say $\pi$ {\it avoids} $\sigma$ whenever $\pi$ contains no subsequence with all of the same pairwise comparisons as $\sigma$.
For example, the permutation 214538769 avoids 312 and 2413, but it has 2586 as a subsequence so it does not avoid 1243.
If $\pi$ avoids $\sigma$ then $\sigma$ is sometimes called a {\it pattern} or a {\it forbidden subsequence} and $\pi$ is sometimes called a {\it restricted permutation} or a {\it pattern-avoiding permutation}.
In this paper we will be interested in permutations which avoid several patterns, so for any set $R$ of permutations we write $S_n(R)$ to denote the elements of $S_n$ which avoid every element of $R$.
For any set $R$ of permutations we take $S_n(R)$ to be the empty set whenever $n < 0$ and we take $S_0(R)$ to be the set containing only the empty permutation.
When $R = \{\pi_1, \pi_2, \ldots, \pi_r\}$ we often write $S_n(R) = S_n(\pi_1, \pi_2, \dots, \pi_r)$.

Pattern avoidance has proved to be a useful language in a variety of seemingly unrelated problems, from singularities of Schubert varieties \cite{LakSong}, to Chebyshev polynomials of the second kind \cite{CW,Krattenthaler,MansourVainshtein2}, to rook polynomials for a rectangular board \cite{MansourVainshtein3}.
There is a particularly fruitful interplay between pattern-avoiding permutations and the study of various sorting algorithms.
One example of this interplay occurs in \cite{W20,W21}, where West studied permutations sortable by two passages through a stack, under the condition that no element of the permutation may be placed on top of a smaller element.
Such permutations are now known as {\em two-stack sortable permutations};  we write ${\cal P}_n$ to denote the set of two-stack sortable permutations in $S_n$.
West characterized two-stack sortable permutations in terms of pattern-avoidance, showing that ${\cal P}_n$ is the set of permutations in $S_n$ which avoid $2341$ and $3241$, except that the latter pattern is allowed whenever it is contained in a subsequence of type $35241$ in the permutation.
West \cite{W20,W21} conjectured that for all $n\geq 0$,
$$|{\cal P}_n|=\frac{2(3n)!}{(n+1)!(2n+1)!}.$$
This conjecture was first proved by Zeilberger in \cite{Z24}.
Zeilberger's proof is analytical, rather than combinatorial.
Later, Dulucq, Gire, and Guibert ~\cite{DG1} gave a combinatorial proof of this conjecture (see also \cite{Br2,BT3,DG7,T18}).

In this paper we consider two-stack sortable permutations which also avoid additional permutations.
Accordingly, for any set $R$ of permutations, we write $\cP_n(R)$ to denote the set of two-stack sortable permutations in $S_n$ which avoid every pattern in $R$, and we write $P_R(x)$ to denote the generating function given by
$$P_R(x) = \sum_{n=0}^\infty |\cP_n(R)| x^n.$$

We will encounter numerous sequences in this paper, many of which will be closely related to two particular sequences:  the sequence of Fibonacci numbers and the sequence of Pell numbers.
We write $F_0, F_1, \ldots$ to denote the sequence of Fibonacci numbers, which are given by $F_0 = 0$, $F_1 = 1$, and $F_n = F_{n-1} + F_{n-2}$ for $n \ge 2$.
We observe that the generating function for the Fibonacci numbers is given by
\begin{equation}
\label{eqn:fibgf}
\sum_{n=0}^\infty F_n x^n = \frac{x}{1 - x - x^2}.
\end{equation}
We also observe that $F_n$ may be interpreted combinatorially as the number of tilings of a $1 \times (n-1)$ rectangle with tiles of size $1 \times 1$ and $1 \times 2$.
We write $p_0, p_1, \ldots$ to denote the sequence of Pell numbers, which are given by $p_0 = 0$, $p_1 = 1$, and $p_n = 2 p_{n-1} + p_{n-2}$ for $n \ge 2$.
We observe that the generating function for the Pell numbers is given by
\begin{equation}
\label{eqn:pellgf}
\sum_{n=0}^\infty p_n x^n = \frac{x}{1 - 2x - x^2}.
\end{equation}
We also observe that $p_n$ may be interpreted combinatorially as the number of tilings of a $1 \times (n-1)$ rectangle with tiles of size $1 \times 1$ and $1 \times 2$, where each $1 \times 1$ tile can be red or blue.
For more information on the Pell numbers, see \cite{FibProblem}, \cite[pp. 122--125]{Beiler}, \cite{Emerson}, and \cite{Horadam}.

In this paper we use generating function techniques to study those
two-stack sortable permutations which avoid (or contain exactly
once) $132$ and which avoid (or contain exactly once) an arbitrary
pattern in $S_k$. In Section \ref{sec2} we describe a method for
enumerating $\cP_n(132,\tau)$ for any $\tau \in \cP_n(132)$. Using
this method, we give several enumerations which involve the
Fibonacci or Pell numbers. In Section \ref{sec3} we use similar
techniques to find generating functions for the set of
permutations in $\cP_n(132)$ which contain another pattern $\tau$
exactly once. Again we give specific examples involving the
Fibonacci or Pell numbers. In Section \ref{sec4} we turn our
attention to the set $\cQ_n$ of permutations in $\cP_n$ which
contain exactly one subsequence of type 132. We give a technique
for finding the generating function for those permutations in
$\cQ_n$ which avoid, or contain exactly once, a permutation
$\tau$. Again we provide specific examples involving the Fibonacci
or Pell numbers. We conclude the paper with examples of other
problems which might be solved using similar techniques.

\section{Two-stack Sortable Permutations Which Avoid $132$ and Another Pattern}
\label{sec2}

In this section we describe a method for enumerating two-stack sortable permutations which avoid 132 and at most one additional pattern and we use our method to enumerate $\cP_n(132,\tau)$ for various $\tau \in \cP_k(132)$.
We observe that since 35241 contains a subsequence of type 132 (namely, 354), for any set $R$ of permutations we have $\cP_n(132, R) = S_n(132, 2341, 3241, R)$.
We begin with an observation concerning the structure of the permutations in $\cP_n(132)$.

\begin{proposition}
\label{prop:prom0}
Fix $n \ge 3$ and suppose $\pi \in \cP_n(132)$.
Then the following hold.
\renewcommand\labelenumi{{\upshape (\roman{enumi}) }}
\begin{enumerate}
\item $\pi^{-1}(n) = 1$, $\pi^{-1}(n) = 2$, or $\pi^{-1}(n) = n$.
\item
The map from $\cP_{n-1}(132)$ to $\cP_n(132)$ given by
$$\pi \mapsto n, \pi$$
is a bijection between $\cP_{n-1}(132)$ and the set of permutations in $\cP_n(132)$ which begin with $n$.
\item
The map from $\cP_{n-2}(132)$ to $\cP_n(132)$ given by
$$\pi \mapsto n-1, n, \pi$$
is a bijection between $\cP_{n-2}(132)$ and the set of permutations in $\cP_n(132)$ whose second entry is $n$.
\item
The map from $\cP_{n-1}(132)$ to $\cP_n(132)$ given by
$$\pi \mapsto \pi,n$$
is a bijection between $\cP_{n-1}(132)$ and the set of permutations in $\cP_n(132)$ which end with $n$.
\end{enumerate}
\end{proposition}
\begin{proof}
(i)
Suppose by way of contradiction that $2 < \pi^{-1}(n) < n$.
Since $\pi$ avoids 132, the elements to the left of $n$ in $\pi$ are all greater than every element to the right of $n$.
Since there are at least two elements to the left of $n$ and at least one element to the right of $n$, there must be a pattern of type 3241 or 2341 in $\pi$ in which $n$ plays the role of the 4.
This contradicts our assumption that $\pi \in \cP_n(132)$.

(ii)
Since the given map is clearly injective it is sufficient to show that $\pi \in \cP_{n-1}(132)$ if and only if $n,\pi \in \cP_n(132)$.
Since $\cP_n(132) = S_n(132, 2341, 3241)$, it is clear that if $n,\pi \in \cP_n(132)$ then $\pi \in \cP_{n-1}(132)$.
To show the converse, suppose $\pi \in \cP_{n-1}(132)$.
If $n,\pi$ contains a pattern of type 132 then the $n$ cannot be involved in the pattern, since it is the largest element of $n,\pi$.
Then $\pi$ contains a pattern of type 132, a contradiction.
By similar arguments for 2341 and 3241 we find that $n,\pi \in \cP_n(132)$, as desired.

(iii),(iv)
These are similar to the proof of (ii).
\end{proof}

\begin{theorem}
\label{thm:Pn132enum}
For all $n \ge 1$,
\begin{equation}
\label{eqn:twostackpell}
|\cP_n(132)| = p_n.
\end{equation}
\end{theorem}
\begin{proof}
For notational convenience, we abbreviate $P(x) = P_{132}(x)$.
By Proposition \ref{prop:prom0}(i), when $n \ge 3$ the elements of $\cP_n(132)$ may be partitioned into three classes:  those which begin with $n$, those whose second entry is $n$, and those whose last entry is $n$.
By Proposition \ref{prop:prom0}(ii), the generating function for those elements which begin with $n$ is $x (P(x) - 1 - x)$.
By Proposition \ref{prop:prom0}(iii), the generating function for those elements whose second entry is $n$ is $x^2 (P(x) - 1)$.
By Proposition \ref{prop:prom0}(iv), the generating function for those elements which end with $n$ is $x (P(x) - 1 - x)$.
Combine these observations to obtain
$$P(x) = 1 + x + 2x^2 + x (P(x) - 1 - x) + x^2 (P(x) -1) + x (P(x) - 1 - x).$$
Solve this equation for $P(x)$ and compare the result with (\ref{eqn:pellgf}) to complete the proof.
\end{proof}

We also give a combinatorial proof of (\ref{eqn:twostackpell}).

\begin{theorem}
For all $n \ge 1$, there exists a constructive bijection between $\cP_n(132)$ and the set of tilings of a $1 \times (n-1)$ rectangle with tiles of size $1 \times 1$ and $1 \times 2$, where each $1 \times 1$ tile can be red or blue.
\end{theorem}
\begin{proof}
Suppose we are given such a tiling;  we construct the corresponding permutation as follows.
Proceed from right to left, placing one number in each box of a tile.
If the rightmost empty tile is a $1 \times 2$ tile then fill it with the two smallest remaining numbers, in increasing order.
If the rightmost empty tile is a red $1 \times 1$ tile then fill it with the largest remaining number.
If the rightmost empty tile is a blue $1 \times 1$ tile then fill it with the smallest remaining number.
When all tiles have been filled, one number will remain.
Place this number in the leftmost position in the permutation.
To obtain a permutation in $\cP_n(132)$, take the (group theoretic) inverse of the permutation constructed by the process above.
It is routine to construct the inverse of the map described above, and thus to verify it is a bijection.
\end{proof}

Though elementary, Proposition \ref{prop:prom0} enables us to easily find $P_{132, \tau}(x)$ for various $\tau$.
For instance, we have the following result involving $\tau = 12\ldots d$.

\begin{theorem}
\label{thm:13212kgf}
\renewcommand\labelenumi{{\upshape (\roman{enumi}) }}
\begin{enumerate}
\item
${\displaystyle P_{132,1}(x) = 1}$.
\item
${\displaystyle P_{132,12}(x) = \frac{1}{1-x}}$.
\item
For all $d \ge 3$,
\begin{equation}
\label{eqn:13212kgf}
P_{132,12\ldots d}(x) = \frac{(1-x) \sum\limits_{r=0}^{d-3} (1-x-x^2)^{r+1} x^{d-3-r} + x^{d-2}}{(1-x)(1-x-x^2)^{d-2}}.
\end{equation}
\end{enumerate}
\end{theorem}
\begin{proof}
(i)
Observe that only the empty permutation avoids 1.

(ii)
Observe that for all $n \ge 0$, the only permutation in $S_n$ which avoids 12 is $n, n-1, \ldots, 2, 1$.

(iii)
We argue by induction on $d$.
By Proposition \ref{prop:prom0}, for all $d \ge 3$ we have
$$P_{132, 12\ldots d}(x) = 1 + x + 2x^2 + x (P_{132, 12\ldots d}(x) - 1 - x) + x^2 (P_{132, 12\ldots d}(x) - 1) + x (P_{132, 12\ldots d-1}(x) - 1 - x).$$
Solve this equation for $P_{132, 12\ldots d}(x)$ to find that for all $d \ge 3$,
\begin{equation}
\label{eqn:13212kpf1}
P_{132, 12\ldots d}(x) = 1 + \frac{x}{1- x- x^2} P_{132, 12\ldots d-1}(x).
\end{equation}
Set $d = 3$ in (\ref{eqn:13212kpf1}), use (ii) to eliminate $P_{132, 12}(x)$, and simplify the result to obtain
$$P_{132, 123}(x) = \frac{(1-x)(1-x-x^2) + x}{(1-x)(1-x-x^2)}.$$
Therefore (\ref{eqn:13212kgf}) holds for $d = 3$.
Moreover, if (\ref{eqn:13212kgf}) holds for $d$ then it is routine using (\ref{eqn:13212kpf1}) to verify that (\ref{eqn:13212kgf}) holds for $d+1$.
\end{proof}

\begin{corollary}
For all $n \ge 1$,
\begin{equation}
\label{eqn:123enum}
|\cP_n(132, 123)| = F_{n+2}-1
\end{equation}
and
\begin{equation}
\label{eqn:1234enum}
|\cP_n(132, 1234)| = 1 - \frac85 F_{n+1} + \frac{n+1}{5} \left( F_{n+3} +
F_{n+1}\right).
\end{equation}
\end{corollary}
\begin{proof}
To prove (\ref{eqn:123enum}), first set $d = 3$ in (\ref{eqn:13212kgf}) and simplify the result to obtain
$$P_{132, 123}(x) = 1 - \frac{1}{1-x} + \frac{x+1}{1-x-x^2}.$$
Compare this last line with (\ref{eqn:fibgf}) to obtain (\ref{eqn:123enum}).
To prove (\ref{eqn:1234enum}), first set $d = 4$ in (\ref{eqn:13212kgf}) and simplify the result to obtain
$$P_{132, 1234}(x) = 1 + \frac{1}{1-x} - \frac{2}{1-x-x^2} + \frac{1}{(1-x-x^2)^2}.$$
It is routine to verify that
$$\frac{1}{(1-x-x^2)^2} = \sum_{n=0}^\infty \left( \left( \frac{1}{5} n + \frac{1}{5} \right) F_n + \left( \frac{3}{5} n + 1 \right) F_{n+1} \right) x^n,$$
and (\ref{eqn:1234enum}) follows.
\end{proof}

Next we consider $P_{132, d12\ldots(d-1)}(x)$.

\begin{theorem}
\label{thm:132k1gf}
\renewcommand\labelenumi{{\upshape (\roman{enumi}) }}
\begin{enumerate}
\item ${\displaystyle P_{132,21}(x) = \frac{1}{1-x}}$.
\item ${\displaystyle P_{132,312}(x) = \frac{(1-x-x^2)(1-x) + x(1+x)}{(1-x)^2}}$.
\item For all $d \ge 4$,
\begin{eqnarray}
\lefteqn{P_{132,d12\ldots(d-1)}(x) =} & & \nonumber \\
\label{eqn:132k1gf}
& &  \frac{(1-x)(1-x-x^2)^{d-2} + (1-x^2) \sum\limits_{r=0}^{d-4} (1-x-x^2)^{r+1} x^{d-3-r} + x^{d-2}(1+x)}{(1-x)^2 (1-x-x^2)^{d-3}}.
\end{eqnarray}
\end{enumerate}
\end{theorem}
\begin{proof}
(i)
Observe that for all $n \ge 0$, the only permutation in $S_n$ which avoids 12 is $1,2,\ldots,n-1,n$.

(ii),(iii)
By Proposition \ref{prop:prom0}, for all $d \ge 3$ we have
\begin{eqnarray*}
\lefteqn{P_{132,d12\ldots(d-1)}(x)} & & \\
&=& 1 + x + 2x^2 + x(P_{132,12\ldots d-1}(x) - 1 - x) + x^2 (P_{132,12\ldots d-1}(x) - 1) \\
& & \hspace{150pt} + x (P_{132,d12\ldots(d-1)}(x) - 1 - x).\\
\end{eqnarray*}
Solve this equation for $P_{132,d12\ldots(d-1)}(x)$ and use Theorem \ref{thm:13212kgf}(ii),(iii) to eliminate the factor $P_{132,12\ldots d-1}(x)$, obtaining (\ref{eqn:132k1gf}).
\end{proof}

\begin{corollary}
We have
\begin{equation}
\label{eqn:312enum}
|\cP_n(132, 312)| = 2n-2 \hspace{30pt} (n \ge 2)
\end{equation}
and
\begin{equation}
\label{eqn:4123enum}
|\cP_n(132, 4123)| = F_{n+4}-2n-2 \hspace{30pt} (n \ge 0).
\end{equation}
\end{corollary}
\begin{proof}
To prove (\ref{eqn:312enum}), first observe from Theorem \ref{thm:132k1gf}(ii) that
$$P_{132,312}(x) = 3 + x + \frac{2}{(1-x)^2} - \frac{4}{1-x}.$$
Now (\ref{eqn:312enum}) follows from the binomial theorem.
The proof of (\ref{eqn:4123enum}) is similar to the proof of (\ref{eqn:312enum}).
\end{proof}

We conclude this section by describing a recursive method for computing $P_{132, \tau}(x)$ for any permutation $\tau \in \cP_k(132)$.
Observe that this allows us to compute $P_{132,\tau}(x)$ for any permutation $\tau \in S_k$, since $P_{132,\tau}(x) = P_{132}(x)$ if $\tau \in S_k$ but $\tau \not\in \cP_k(132)$.

\begin{theorem}
\label{thm:132tau}
Fix $k \ge 3$ and $\tau \in \cP_k(132)$.
Observe by Proposition \ref{prop:prom0} that exactly one of the following holds.
\begin{enumerate}
\item[(a)]
There exists $\tau' \in \cP_{k-1}(132)$ such that $\tau = k,\tau'$.
\item[(b)]
There exists $\tau' \in \cP_{k-2}(132)$ such that $\tau = k-1,k,\tau'$.
\item[(c)]
There exists $\tau' \in \cP_{k-1}(132)$ such that $\tau = \tau',k$.
\end{enumerate}
Then the following also hold.
\renewcommand\labelenumi{{\upshape (\roman{enumi}) }}
\begin{enumerate}
\item
If (a) holds then
\begin{equation}
P_{132,\tau}(x) = \frac{1-x-x^2 + x (1+x) P_{132,\tau'}(x)}{1-x}.
\end{equation}
\item
If (b) holds then
\begin{equation}
P_{132,\tau}(x) = \frac{1-x-x^2 + x^2 P_{132,\tau'}(x)}{1-2x}.
\end{equation}
\item
If (c) holds then
\begin{equation}
P_{132,\tau}(x) = 1 + \frac{x P_{132,\tau'}(x)}{1-x-x^2}.
\end{equation}
\end{enumerate}
\end{theorem}
\begin{proof}
The proofs of (i)--(iii) are similar to the proofs of Theorem \ref{thm:13212kgf}(iii) and \ref{thm:132k1gf}(iii).
\end{proof}

\begin{corollary}
\begin{equation}
\label{eqn:3412enum}
|\cP_n(132,3412)| = 3 \cdot 2^{n-2} - 1 \hspace{30pt} (n \ge 2)
\end{equation}
\begin{equation}
\label{eqn:45123enum}
|\cP_n(132,45123)| = 3 \cdot 2^{n-1} - F_{n+3} + 1 \hspace{30pt} (n \ge 1)
\end{equation}
\begin{equation}
\label{eqn:561234enum}
|\cP_n(132,561234)| = 3 \cdot 2^n - 1 - \frac{6}{5} F_{n+1} - \frac{n+1}{5} (F_{n+4}+F_{n+2}) \hspace{30pt} (n \ge 1)
\end{equation}
\end{corollary}
\begin{proof}
To prove (\ref{eqn:3412enum}), first use Theorems \ref{thm:132tau}(ii) and \ref{thm:13212kgf}(ii) to find that
$$P_{132, 3412}(x) = \frac{5}{4} + \frac{1}{2} x - \frac{1}{1-x} + \frac{3}{4} \frac{1}{1-2x}.$$
Now (\ref{eqn:3412enum}) is immediate.
The proofs of (\ref{eqn:45123enum}) and (\ref{eqn:561234enum}) are similar to the proof of (\ref{eqn:3412enum}).
\end{proof}

\begin{proposition}
For all $d\ge 2$,
\begin{equation}
\label{eqn:132k21gf}
P_{132, d\ldots 21}(x) = \frac{(1-x-x^2) \sum\limits_{i=0}^{d-3} x^i (1 + x)^i (1-x)^{d-i-2} + x^{d-2} (1+x)^{d-2}}{(1-x)^{d-1}}.
\end{equation}
\end{proposition}
\begin{proof}
We argue by induction on $d$.
When $d = 2$, line (\ref{eqn:132k21gf}) reduces to Theorem \ref{thm:132k1gf}(ii).
If (\ref{eqn:132k21gf}) holds for a given $d \ge 2$ then it is routine using Theorem \ref{thm:132tau}(i) to show that (\ref{eqn:132k21gf}) holds for $d + 1$.
\end{proof}

\begin{corollary}
We have
\begin{equation}
\label{eqn:321enum}
|\cP_n(132, 321)| = 2n-2 \hspace{30pt} (n \ge 2),
\end{equation}
\begin{equation}
\label{eqn:4321enum}
|\cP_n(132, 4321)| = 2 n^2 - 8n + 11 \hspace{30pt} (n \ge 3),
\end{equation}
and
\begin{equation}
\label{eqn:54321enum}
|\cP_n(132, 54321)| = \frac{4}{3} n^3 - 12 n^2 + \frac{128}{3} n - 52 \hspace{30pt} (n \ge 4).
\end{equation}
\end{corollary}
\begin{proof}
To prove (\ref{eqn:321enum}), first set $d = 3$ in (\ref{eqn:132k21gf}) and simplify the result to find
$$P_{132, 321}(x) = 3 + x - \frac{4}{1-x} + \frac{2}{(1-x)^2}.$$
Now (\ref{eqn:321enum}) follows from the binomial theorem.
The proofs of (\ref{eqn:4321enum}) and (\ref{eqn:54321enum}) are similar to the proof of (\ref{eqn:321enum}).
\end{proof}

\section{Two-stack Sortable Permutations Which Avoid $132$ and Contain Another Pattern}
\label{sec3}

Fix $d \ge 1$ and set
$$\cP(132) = \bigcup_{n\ge 0} \cP_n(132).$$
Inspired by results such as those found in \cite{BCS} and \cite{MS}, we consider in this section the generating function for permutations in $\cP(132)$ according to the number of patterns of type $12\ldots d$, $d\ 1 2 \ldots d-1$, or $d\ d-1 \ldots 2 1$ they contain.
We begin by setting some notation.

\begin{definition}
For any permutation $\tau \in \cP_k(132)$ and any $r \ge 1$, let $b_{n,r}$ denote the number of permutations in $\cP_n(132)$ which contain exactly $r$ subsequences of type $\tau$.
Then we write
\begin{equation}
\label{eqn:Bntau}
B_\tau^r(x) = \sum_{n=0}^\infty b_{n,r} x^n.
\end{equation}
\end{definition}

We now consider the case in which $\tau = 12 \ldots d$.

\begin{definition}
For any permutation $\pi$ and any $d \ge 1$, we write $12\ldots d(\pi)$ to denote the number of subsequences of type $12\ldots d$ in $\pi$ and we write $|\pi|$ to denote the length of $\pi$.
\end{definition}

\begin{theorem}
\label{thm:12dgf}
Let $x_1, x_2, \ldots$ denote indeterminates.
Then we have
\begin{equation}
\label{eqn:12dgf}
\sum_{\pi\in\mathcal{P}(132)}\prod\limits_{d\geq1}x_d^{12\dots
d(\pi)}=1+\sum_{n\geq1}\frac{\prod\limits_{j\geq1}x_j^{{{n}\choose{j}}}}
{\prod\limits_{m=1}^n\left(1-\prod\limits_{j\geq1}x_j^{{{m-1}\choose{j-1}}}-\prod\limits_{j\geq1}x_j^{2{{m-1}\choose{j-1}}}x_{j+1}^{{{m-1}\choose{j-1}}}\right)}.
\end{equation}
\end{theorem}
\begin{proof}
For notational convenience, set
$$A(x_1, x_2, \ldots) = \sum_{\pi \in \cP(132)} \prod_{d \ge 1} x_d^{12\ldots d(\pi)}.$$
By Proposition \ref{prop:prom0},
\begin{eqnarray*}
A(x_1, x_2, \ldots) &=& 1 + x_1 + x_1^2 x_2 + x_1^2 + x_1 \left( A(x_1, x_2, \ldots) - 1 - x_1 \right)\\[2ex]
& & {\ }+ x_1^2 x_2 \left( A(x_1, x_2, \ldots ) - 1 \right) + x_1 \left( A(x_1 x_2, x_2 x_3, \ldots)-1 -x_1 x_2\right) \\
\end{eqnarray*}
Solve this equation for $A(x_1, x_2, \ldots)$ to obtain
$$A(x_1, x_2, \ldots) = 1 + \frac{x_1 A(x_1 x_2, x_2 x_3, \ldots)}{1 - x_1 - x_1^2 x_2}.$$
Iterate this recurrence relation to obtain (\ref{eqn:12dgf}).
\end{proof}

For a given $d \ge 1$, we can use (\ref{eqn:12dgf}) to obtain the generating function for $\cP(132)$ according to length and number of subsequences of type $12\ldots d$.

\begin{proposition}
For all $d \ge 3$,
\begin{equation}
\label{eqn:B11dgf}
B_{12\ldots d}^1(x) = \frac{x^d}{(1-x)^2 (1-x-x^2)^{d-2}}.
\end{equation}
\end{proposition}
\begin{proof}
In (\ref{eqn:12dgf}), set $x_1 = x$, $x_d = y$, and $x_i = 1$ for $i \neq 1, d$.
Expand the resulting expressing in powers of $y$;  we wish to find the coefficient of $y$.
Observe that only the terms in the sum for which $n = d-1$ or $n = d$ contribute to this coefficient.
These terms are
$$\frac{x^{d-1}}{(1-x-x^2)^{d-2} (1-x-x^2y)}$$
and
$$\frac{x^d y}{(1-x-x^2)^{d-2} (1-x-x^2 y) (1-xy-x^2 y^2)}$$
respectively, and (\ref{eqn:B11dgf}) follows.
\end{proof}

\begin{corollary}
For all $n \ge 0$, the number of permutations in $\cP_n(132)$ which contain exactly one subsequence of type 123 is $F_{n+2} - n - 1$.
\end{corollary}
\begin{proof}
Set $d = 3$ in (\ref{eqn:B11dgf}) and simplify the result to find
$$B_{123}^1(x) = \frac{x+1}{1-x-x^2} - \frac{1}{(1-x)^2}.$$
Now the result follows from (\ref{eqn:fibgf}) and the binomial theorem.
\end{proof}

As another application of (\ref{eqn:12dgf}), we now find the generating function for $\cP(132)$ according to length and number of right to left maxima.
To do this, we first find the number of right to left maxima in a given permutation $\pi$ in terms of $12\ldots d(\pi)$.

\begin{definition}
For any $\pi \in S_n$, we write $rmax(\pi)$ to denote the number of right to left maxima in $\pi$.
\end{definition}

\begin{proposition}
For all $\pi \in S_n(132)$, we have
\begin{equation}
\label{eqn:rmax12d}
rmax(\pi) = \sum_{d = 1}^n (-1)^{d+1} 12\ldots d(\pi).
\end{equation}
\end{proposition}
\begin{proof}
Fix $\pi \in S_n(132)$ and fix $i$, $1 \le i \le n$.
We consider the contribution of those increasing subsequences of $\pi$ which begin at $\pi(i)$ to the sum on the right side of (\ref{eqn:rmax12d}).
If $\pi(i)$ is a right to left maxima, then the only increasing subsequence which begins at $\pi(i)$ has length one.
Therefore each right to left maxima contributes one to the sum.
If $\pi(i)$ is not a right to left maxima then we observe that because $\pi$ avoids 132, the elements to the right of $\pi(i)$ which are larger than $\pi(i)$ are in increasing order.
Therefore the contribution of $\pi(i)$ to the sum is
$$\sum_{i=0}^k (-1)^i {{k} \choose {i}} = 0,$$
where $k$ is the number of elements in $\pi$ larger than $\pi(i)$ to the right of $\pi(i)$.
Combine these observations to obtain (\ref{eqn:rmax12d}).
\end{proof}

\begin{proposition}
We have
\begin{equation}
\label{eqn:lengthrmaxgf}
\sum_{\pi \in \cP(132)} x^{|\pi|} y^{rmax(\pi)} = 1 + \frac{xy (1-x-x^2)}{(1-xy-x^2y)(1-2x-x^2)}.
\end{equation}
\end{proposition}
\begin{proof}
In (\ref{eqn:12dgf}), set $x_1 = xy$ and $x_i = y^{(-1)^{i+1}}$ for $i \ge 2$.
Use the facts
$$\sum_{i=0}^n (-1)^i {{n}\choose{i}} = \delta_{n1}$$
and
$$\sum_{n=1}^\infty \frac{x^n}{(1-x-x^2)^{n-1}} = \frac{x (1 -x -x^2)}{1-2x-x^2}$$
to simplify the result.
\end{proof}

\begin{corollary}
For all $r \ge 1$, let $a_{r,n}$ denote the number of permutations in $\cP_n(132)$ with exactly $r$ right to left maxima.
Then
\begin{equation}
\label{eqn:arngf}
\sum_{n=0}^\infty a_{r,n} x^n = \frac{x^r (1+x)^{r-1} (1-x-x^2)}{1-2x-x^2}.
\end{equation}
In particular,
\begin{equation}
\label{eqn:a1n}
a_{1,n} = p_{n-1} \hspace{30pt} (n \ge 2),
\end{equation}
\begin{equation}
\label{eqn:a2n}
a_{2,n} = p_{n-2} + p_{n-3} \hspace{30pt} (n \ge 4),
\end{equation}
and
\begin{equation}
\label{eqn:a3n}
a_{3,n} = 2 p_{n-3} \hspace{30pt} (n \ge 6).
\end{equation}
\end{corollary}
\begin{proof}
To obtain (\ref{eqn:arngf}), first observe that
$$\frac{1}{1-xy-x^2y} = \sum_{n=0}^\infty x^n (x+1)^n y^n.$$
Combine this with (\ref{eqn:lengthrmaxgf}) to find that
$$\sum_{\pi \in \cP(132)} x^{|\pi|} y^{rmax(\pi)} = \frac{x(1-x-x^2)}{1-2x-x^2} \sum_{n=0}^\infty x^n y^{n+1} (1+x)^n.$$
Take the coefficient of $y^r$ in this last line to obtain (\ref{eqn:arngf}).
To obtain (\ref{eqn:a1n}), set $r = 1$ in (\ref{eqn:arngf}) and compare the result with (\ref{eqn:pellgf}).
The proofs of (\ref{eqn:a2n}) and (\ref{eqn:a3n}) are similar to the proof of (\ref{eqn:a1n}).
\end{proof}

We remark that (\ref{eqn:a1n}) and (\ref{eqn:a2n}) can also be obtained directly from Proposition \ref{prop:prom0} and Theorem \ref{thm:Pn132enum}.
For instance, if $\pi \in \cP_n(132)$ has exactly one right to left maxima then it must end in $n$.
By Proposition \ref{prop:prom0}(iii) and Theorem \ref{thm:Pn132enum}, there are exactly $p_{n-1}$ such permutations in $\cP_n(132)$, and (\ref{eqn:a1n}) follows.
A similar but slightly more involved argument proves (\ref{eqn:a2n}).

Next we find $B_{d 1 2 \ldots d-1}^1(x)$.

\begin{theorem}
\label{thm:B1d12}
\renewcommand\labelenumi{{\upshape (\roman{enumi}) }}
\begin{enumerate}
\item
\begin{equation}
\label{eqn:B1312}
B_{312}^1(x) = \frac{x^3}{(1-x)^2}
\end{equation}
\item
For all $d \ge 4$,
\begin{equation}
\label{eqn:B1d12}
B_{d 1 2 \ldots d-1}^1(x) = \frac{x^d}{(1-x)^3 (1 - x -x^2)^{d-4}}.
\end{equation}
\end{enumerate}
\end{theorem}
\begin{proof}
(i)
Fix $\pi \in \cP_n(132)$ such that $\pi$ contains exactly one subsequence of type 312.
We consider the three cases of Proposition \ref{prop:prom0}.
If $\pi(1) = n$ then $n \ge 3$ and $\pi = n, n-2, n-1, n-3, n-4, \ldots, 2, 1$.
If $\pi(2) = n$ and $\pi(1) = n-1$ then $\pi$ does not contain exactly one subsequence of type 312.
If $\pi = \pi', n$ then $\pi' \in \cP_{n-1}(132)$ and $\pi'$ contains exactly one subsequence of type 312.
Combine these observations to find
$$B^1_{312}(x) = \frac{x^3}{(1-x)} + x B^1_{312}(x).$$
Solve this equation for $B_{312}^1(x)$ to obtain (\ref{eqn:B1312}).

(ii)
Fix $\pi \in \cP_n(132)$ such that $\pi$ contains exactly one subsequence of type $d 1 2 \ldots d-1$.
We consider the three cases of Proposition \ref{prop:prom0}.
If $\pi(1) = n$ then $\pi(n) = n-1$, since $\pi$ avoids 132, 2341, and contains exactly one subsequence of type $d 1 2 \ldots d-1$.
Moreover, it is routine to show that the map $\pi' \mapsto n, \pi', n-1$ is a bijection between those permutations in $\cP_{n-2}(132)$ which contain exactly one subsequence of type $12\ldots d-2$ and those permutations in $\cP_n(132)$ which contain exactly one subsequence of type $d 1 2 \ldots d-1$.
If $\pi(2) = n$ and $\pi(1) = n-1$ then $\pi$ does not contain exactly one subsequence of type $d 1 2 \ldots d-1$.
If $\pi(n) = n$ then $\pi = \pi', n$, where $\pi' \in \cP_{n-1}(132)$ and $\pi'$ contains exactly one subsequence of type $d 1 2 \ldots d-1$.
Combine these observations to find
$$B_{d12\ldots d-1}^1(x) = x^2 B^1_{12\ldots d-2}(x) + x B^1_{d12\ldots d-1}(x).$$
Now use (\ref{eqn:B11dgf}) to eliminate $B_{12\ldots d-2}^1(x)$ and solve for $B^1_{d 1 2 \ldots d-1}(x)$ to obtain (\ref{eqn:B1d12}).
\end{proof}

\begin{corollary}
For all $n \ge 1$, the number of permutations in $\cP_n(132)$ which contain exactly one subsequence of type 51234 is
$$F_{n+2} - {{n+1} \choose {2}}.$$
\end{corollary}
\begin{proof}
Set $d = 5$ in (\ref{eqn:B1d12}) and simplify the result to find
$$B^1_{51234}(x) = \frac{x+1}{1 - x - x^2} - \frac{3}{1-x} + \frac{2}{(1-x)^2} - \frac{1}{(1-x)^3}.$$
Now the result follows from (\ref{eqn:fibgf}) and the binomial theorem.
\end{proof}

We conclude this section by finding $B_{d\ d-1\ \ldots 21}^1(x)$.

\begin{theorem}
For all $d \ge 1$,
\begin{equation}
\label{eqn:B1d21}
B^1_{d\ d-1\ \ldots 2 1}(x) = \frac{x^d}{(1-x)^{d-1}}.
\end{equation}
\end{theorem}
\begin{proof}
This is similar to the proof of Theorem \ref{thm:B1d12}.
\end{proof}

\section{Two-stack Sortable Permutations Which Contain 132 Exactly Once}
\label{sec4}

We now turn our attention to another set of permutations in $\cP_n$.

\begin{definition}
For all $n \ge 0$, we write $\cQ_n$ to denote the set of permutations in $\cP_n$ which contain exactly one subsequence of type 132.
\end{definition}

Using the methods of the previous section, one can obtain generating functions for those permutations in $\cQ_n$ which avoid, or contain exactly once, any permutation $\tau \in S_k$.
In this section we illustrate these derivations with examples in which the resulting generating functions are given in terms of the generating functions for the Pell or Fibonacci numbers.
We begin with an analogue of Proposition \ref{prop:prom0}.

\begin{proposition}
\label{prop:prom1}
Fix $n \ge 3$ and suppose $\pi \in \cQ_n$.
Then the following hold.
\renewcommand\labelenumi{{\upshape (\roman{enumi}) }}
\begin{enumerate}
\item
$\pi(1) = n$, $\pi(2) = n$, $\pi(n) = n$, or $\pi = 3142$.
\item
The map from $\cQ_{n-1}$ to $\cQ_n$ given by
$$\pi \mapsto n, \pi$$
is a bijection between $\cQ_{n-1}$ and the set of permutations in $\cQ_n$ which begin with $n$.
\item
The map from $\cQ_{n-2}$ to $\cQ_n$ given by
$$\pi \mapsto n-1, n, \pi$$
is a bijection between $\cQ_{n-2}$ and the set of permutations in $\cQ_n$ whose second entry is $n$ and in which $n$ does not take part in the subsequence of type 132.
\item
The map from $\cQ_{n-1}$ to $\cQ_n$ given by
$$\pi \mapsto \pi,n$$
is a bijection between $\cQ_{n-1}$ and the set of permutations in $\cQ_n$ which end with $n$.
\item
The map from $\cP_{n-2}$ to $\cQ_n$ given by
$$\pi \mapsto n-2, n, \pi',$$
in which $\pi'$ is the permutation of $1,2,\ldots, n-3, n-1$ which results from replacing $n-2$ with $n-1$ in $\pi$, is a bijection between $\cP_{n-2}$ and the set of permutations in $\cQ_n$ whose second entry is $n$ and in which $n$ takes part in the subsequence of type 132.
\end{enumerate}
\end{proposition}
\begin{proof}
(i)
Suppose $2 < \pi^{-1}(n) < n$.
We consider two cases:  $n$ takes part in the subsequence of type 132 or $n$ does not take part in the subsequence of type 132.

If $n$ does not take part in the subsequence of type 132 then the elements to the left of $n$ are all greater than every element to the right of $n$.
Since there are at least two elements to the left of $n$ and at least one element to the right of $n$, there must be a pattern of type 3241 or 2341 in which $n$ plays the role of the 4.
This contradicts our assumption that $\pi \in \cQ_n \subseteq \cP_n.$

Now suppose the subsequence of type 132 in $\pi$ is $a, n, b$.
Since $\pi$ contains no other subsequence of type 132, all elements to the left of $n$ other than $a$ are greater than $b$ and all elements to the right of $n$ other than $b$ are less than $a$.
Observe that if there are additional elements both to the right and left of $n$ then one of these elements on each side of $n$, together with $a$ and $n$, form a pattern of type 2341 or 3241.
Therefore, there can only be additional elements on one side of $n$.
Since we have assumed $2 < \pi^{-1}(n) < n$, there must be at least one additional element to the left of $n$ and no additional elements to the right of $n$.
If there are two (or more) additional elements to the left of $n$ then they combine with $n$ and $b$ to form a 3241 or a 2341 pattern, so there must be exactly one additional element to the left of $n$.
Now the only possibilities are $\pi = 3142$ and $\pi = 1342$.
But 1342 has more than one subsequence of type 132, so we must have $\pi = 3142$, as desired.

(ii)--(v)
These are similar to the proof of Proposition \ref{prop:prom0}(ii).
\end{proof}

Using Proposition \ref{prop:prom1}, we now find the generating function for $\cQ_n$.

\begin{theorem}
We have
\begin{equation}
\label{eqn:Qngf}
\sum_{n=0}^\infty |\cQ_n| x_n = \frac{x^3 (1 + x - 2x^2 - x^3)}{(1 -2x - x^2)^2}.
\end{equation}
Moreover,
\begin{equation}
\label{eqn:Qnenum}
|\cQ_n| = \frac{1}{4}\left( 51 n p_n - 145 p_n - 21n p_{n+1} + 60 p_{n+1}\right) \hspace{30pt} (n \ge 3).
\end{equation}
\end{theorem}
\begin{proof}
For notational convenience, set
$$Q(x) = \sum_{n=0}^\infty |\cQ_n| x^n.$$
To obtain (\ref{eqn:Qngf}), observe that by Proposition \ref{prop:prom1} we have
$$Q(x) = 2x Q(x) + x^2 Q(x) + x^2 P_{132}(x) + x^4.$$
Use (\ref{eqn:twostackpell}) and (\ref{eqn:pellgf}) to eliminate $P_{132}(x)$ and solve the resulting equation for $Q(x)$ to obtain (\ref{eqn:Qngf}).

To obtain (\ref{eqn:Qnenum}), first observe that
$$Q(x) = -x^2 + 4x - 15 + \frac{87x - 36}{(1-2x-x^2)^2} - \frac{49x - 51}{1-2x-x^2}.$$
Now observe that
$$\frac{1}{(1-2x-x^2)^2} = \sum_{n=0}^\infty \frac{1}{4} \left( (n+1)p_n + (3n+4)p_{n+1}\right) x^n.$$
Combine these observations with (\ref{eqn:pellgf}) to obtain (\ref{eqn:Qnenum}).
\end{proof}

Proposition \ref{prop:prom1} enables us to find the generating function for those permutations in $\cQ_n$ which avoid, or contain exactly once, any permutation $\tau$.
We illustrate how this is done by considering those permutations in which $\tau = 12\ldots d$ appears exactly once.
This case is particularly interesting, since the result generating function is expressed in terms of the generating function for the Fibonacci numbers.

\begin{definition}
For any permutation $\tau \in S_k$, let $c_n$ denote the number of permutations in $\cQ_n$ which contain exactly one subsequence of type $\tau$.
We write $D^1_\tau(x)$ to denote the generating function given by
$$D^1_\tau(x) = \sum_{n=0}^\infty c_n x^n.$$
\end{definition}

\begin{theorem}
For all $d \ge 2$,
\begin{equation}
\label{eqn:D11dgf}
D^1_{12\ldots d}(x) = \frac{(d-2)x^{d+2}}{(1-x)^2 (1-x-x^2)^{d-1}}.
\end{equation}
\end{theorem}
\begin{proof}
The case $d = 2$ is immediate, so we assume $d > 2$ and argue by induction on $d$.
By Proposition \ref{prop:prom1}, we find
$$D^1_{12\ldots d}(x) = x D^1_{12\ldots d}(x) + x^2 D^1_{12\ldots d}(x) + x D^1_{12\ldots d-1}(x) + x^2 B^1_{12\ldots d}(x).$$
Use (\ref{eqn:B11dgf}) to eliminate $B^1_{12\ldots d}(x)$ and solve the resulting equation for $D^1_{12\ldots d}(x)$, obtaining
$$D^1_{12\ldots d}(x) = \frac{1}{1-x-x^2} \left( x D^1_{12\ldots d}(x) + \frac{x^{d+2}}{(1-x)^2 (1-x-x^2)^{d-2}}\right).$$
Use induction to eliminate $D^1_{12\ldots d-1}(x)$ and (\ref{eqn:D11dgf}) follows.
\end{proof}

\begin{corollary}
The number of permutations in $\cQ_n$ which contain exactly one subsequence of type 123 is
$$\frac{n}{5} \left( F_{n+1} + F_{n-1} \right) - \frac{2}{5} \left( 5 F_{n+1} - F_{n-1} \right) +n + 2 \hspace{30pt} (n \ge 4).$$
\end{corollary}
\begin{proof}
This is similar to the proof of (\ref{eqn:1234enum}).
\end{proof}

\section{Directions for Future Work}

In this section we present several directions in which this work may be generalized, using similar techniques.

\begin{enumerate}
\item
In Section \ref{sec2} we enumerated two-stack sortable permutations which avoid 132 and one additional pattern.
Using the same techniques, one ought to be able to enumerate two-stack sortable permutations which avoid 132 and two or more additional patterns.
For instance, for any $d \ge 2$, it should be possible using our techniques to find the generating function $P_{132,123\dots d, 213\dots
d}(x)$.
Moreover, we expect that this generating function will be given in terms of the generating function for the Fibonacci numbers.

\item
In Section \ref{sec4} we enumerated two-stack sortable permutations which contain exactly one subsequence of type 132 and avoid (or contain exactly once) another pattern.
One ought to be able to use similar techniques to enumerate permutations which contain exactly $r$ subsequences of type 132, for a given $r$.

\item
In Section \ref{sec3} we found generating functions for two-stack sortable permutations with respect to various statistics on permutations.
Using similar techniques, one ought to be able to find generating functions for two-stack sortable permutations with respect to additional statistics, such as rises, descents, right to left maxima, right to left minima, left to right maxima, and left to right minima.
\end{enumerate}

\end{document}